\documentstyle[12pt,amscd,amssymb]{amsart}
\newtheorem{thm}{Theorem}[section]
\newtheorem{lemma}[thm]{Lemma}

\newtheorem{prop}[thm]{Proposition}

\newtheorem{cor}[thm]{Corollary}
\newtheorem{defin}[thm]{Definition}
\begin{document}

\newcommand{\C}{{\mathbb C}}
\newcommand{\I}{{\o I}}
\newcommand{\CH}{{\mathbb H}_{\C}^2}
\renewcommand{\P}{{\mathbb P}}
\newcommand{\R}{{\mathbb R}}
\newcommand{\Z}{{\mathbb Z}}
\newcommand{\Os}{{\mathcal O}}
\newcommand{\zz}{{\mathcal Z}}
\newcommand{\z}{{\zz}^1(M)}
\newcommand{\M}{{\mathcal M}}
\newcommand{\la}{\langle}
\newcommand{\ra}{\rangle}
\newcommand{\h}{{\mathcal H}^3_{\R}}
\newcommand{\hH}{{\mathscr H}}
\newcommand{\hi}{h^{-1}}
\renewcommand{\o}{\operatorname}
\renewcommand{\Im}{\o{Im}}
\renewcommand{\Re}{\o{Re}}
\newcommand{\End}{\o{End(E)}}
\newcommand{\Hol}{\o{Hol}(E)}
\newcommand{\Her}{\o{Her}(E)}
\newcommand{\fc}{{\mathcal F}(E)} 
\newcommand{\fcu}{{\mathcal F}_u(E)} 
\newcommand{\fcl}{{\mathcal F}_l(E)} 
\newcommand{\ga}{{\mathcal G}(E)} 
\newcommand{\gau}{{\mathcal G}(E)} 
\newcommand{\gal}{{\mathcal G}_{\C}(E)} 
\newcommand{\Hig}{\o{Higgs}(V)}
\newcommand{\Ext}{\o{Ext}}
\newcommand{\Hom}{\o{Hom}}
\newcommand{\Up}{\o{U}(p)}
\newcommand{\up}{{\frak u}(p)}
\newcommand{\PUp}{\o{PU}(p)}
\newcommand{\Uq}{\o{U}(q)}
\newcommand{\uq}{{\frak u}(q)}
\newcommand{\PUq}{\o{PU}(q)}
\newcommand{\Un}{\o{U}(n)}
\newcommand{\un}{{\frak u}(n)}
\newcommand{\PUn}{\o{PU}(n)}
\newcommand{\Uoo}{\o{U}(1,1)}
\newcommand{\uoo}{{\frak u}(1,1)}
\newcommand{\PUoo}{\o{PU}(1,1)}
\newcommand{\Ut}{\o{U}(2)}
\newcommand{\ut}{{\frak u}(2)}
\newcommand{\Uto}{\o{U}(2,1)}
\newcommand{\PUto}{\o{PU}(2,1)}
\newcommand{\uto}{{\frak u}(2,1)}
\newcommand{\Upq}{\o{U}(p,q)}
\newcommand{\Upo}{\o{U}(p,1)}
\newcommand{\upq}{{\frak u}(p,q)}
\newcommand{\upo}{{\frak u}(p,1)}
\newcommand{\PUpq}{\o{PU}(p,q)}
\newcommand{\PUpo}{\o{PU}(p,1)}
\newcommand{\Uo}{\o{U}(1)}
\newcommand{\uo}{{\frak u}(1)}
\newcommand{\GLn}{\o{GL}(n,\C)}
\newcommand{\gln}{{\frak gl}(n,\C)}
\newcommand{\PGLn}{\o{PGL}(n,\C)}
\newcommand{\PSLt}{\o{PSL}(3,\C)}
\newcommand{\rk}{\o{rank}}
\newcommand{\Map}{\o{Map}}
\newcommand{\db}{\overline{\partial}}
\renewcommand{\d}{\partial}
\newcommand\ka{K\"ahler~}
\newcommand\tM{\tilde{M}}
\newcommand\ca{{\mathcal A}}
\newcommand\cb{{\mathcal B}}
\newcommand\cc{{\mathcal C}}
\renewcommand{\i}{{\o i}}
\renewcommand{\j}{{\o j}}
\renewcommand{\k}{{\o k}}
\newcommand\Iso{{\o{Iso}}}
\newcommand\Obj{{\o{Obj}}}
\newcommand\hk{{hyperk\"ahler~}}
\newcommand{\E}{{\cal E}}
\newcommand{\F}{{\cal F}}
\newcommand{\G}{{\cal G}}
\renewcommand{\H}{{\o H}}
\newcommand{\K}{{\cal K}}
\renewcommand{\L}{{\cal L}}
\newcommand{\N}{{\cal N}}
\newcommand{\U}{{\cal U}}
\newcommand{\V}{{\cal V}}
\newcommand{\W}{{\cal W}}
\renewcommand{\O}{{\cal O}}
\newcommand{\Q}{{\cal Q}}

\title
{The Moduli of Flat U($p,1$) Structures on Riemann Surfaces}
\author{Eugene Z. Xia}
\address{ 
Department of Mathematics,
University of Massachusetts, Amherst, MA 01003-4515}
\email{xia@@math.umass.edu}
\date{\today} 
\subjclass{
14D20 (Algebraic Moduli Problems, Moduli of Vector Bundles),
14H60 (Vector Bundles on Curves)}
\keywords{
Algebraic curves, Moduli schemes, Higgs bundles}
\begin{abstract}
For a compact Riemann surface $X$ of genus $g > 1$, 
$\Hom(\pi_1(X), \Upo)/\Upo$ is the moduli space of 
flat $\Upo$-connections on $X$.
There is an integer invariant, $\tau$, the Toledo invariant
associated with each element in $\Hom(\pi_1(X), \Upo)/\Upo$. 
If $q = 1$, then $-2(g-1) \le \tau \le 2(g-1)$.
This paper shows that 
$\Hom(\pi_1(X), \Upo)/\Upo$ has one connected component
corresponding to each $\tau \in 2\Z$ with
$-2(g-1) \le \tau \le 2(g-1)$.
Therefore the total number of connected components is $2(g-1) + 1$.
\end{abstract}
\maketitle


\section{Introduction and Results}
Let $X$ be a smooth projective curve over $\C$ with genus $g>1$.
The deformation space 
$$
\C {\cal M}_B = \Hom^+(\pi_1(X), \GLn)/\GLn
$$
is the space of equivalence classes of 
semi-simple $\GLn$-representations of the fundamental
group $\pi_1(X)$.
This is the $\GLn$-Betti moduli space on $X$. 
A theorem of Corlette, Donaldson, Hitchin and Simpson 
relates $\C {\cal M}_B$ to
two other moduli spaces---the 
$\GLn$-de Rham and the $\GLn$-Dolbeault moduli spaces, respectively
\cite{Co1,Do1, Hi1, Si1}.
The Dolbeault moduli space consists of holomorphic objects (Higgs
bundles) over $X$; therefore, the classical results of
analytic and algebraic geometry can
be applied to the study of the Dolbeault moduli space.

Since $\Upo \subset \GLn$, $\C {\cal M}_B$ contains the space
$$
{\cal M}_B = \Hom^+(\pi_1(X), \Upo)/\Upo.
$$
The space ${\cal M}_B$ will be referred to as the $\Upo$-Betti moduli space.

The Betti moduli spaces are of great interest in
geometric topology and uniformization.
When $p=q=1$,
Goldman analyzed ${\cal M}_B$
and determined the number of its connected
components \cite{Go1}.
Hitchin subsequently considered this moduli space
from the Higgs bundle perspective
and determined its topology \cite{Hi1}.
The case of $p=2,q=1$ was treated in \cite{Xi3}.
Other related results have been obtained in \cite{Go2, Xi1, Xi2}.
In this paper, we treat the general case of 
$q=1$ and determine its number of connected components.

Each element in ${\cal M}_B$ is associated with a Toledo invariant
$\tau \in 2 \Z$ \cite{Do1,To1,To2}.
The main result presented here is the following:
\begin{thm} \label{thm:main}
$\Hom^+(\pi_1(X), \Upo)/\Upo$ has one connected component
for each $\tau \in 2\Z$ with 
$-2(g-1) \le \tau \le 2(g-1)$.  Therefore the total number of
connected components is $2(g-1) + 1$.
\end{thm}
\subsection{The $\PUpo$-representations}  
One hopes
to prove a similar result for the case of $\PUpo$-representations i.e.
there is one connected component for each possible $\tau$.
Here the Toledo invariant is also bounded as
$-2(g-1) \le \tau \le 2(g-1)$, but takes on values in $\frac{2}{p+1} \Z$
\cite{Xi3}.  Again the cases of $p = 1,2$ have been treated in
\cite{Go1, Hi1, Xi3}.
At present, the case of general $p$ seems rather difficult.  
The reason we are able to treat
the $\Upo$-representations is due to the presence of 
certain reducible $\Upo$-representations.
These representations correspond to semi-stable
Higgs bundles that can be constructed rather explicitly.

\section{Backgrounds and Preliminaries}
The coarse moduli space $M_{r,d}$ of 
semi-stable vector bundles on $X$ of rank $r$ and degree $d$
exists and has dimension $r^2(g-1) +1$ \cite{Se1}.  

\subsection{The $\Upo$-Higgs bundles} 
Each element $\rho \in \Hom(\pi_1(X), \GLn)$ acts on $\C^n$
via the standard representation of $\GLn$.  The representation $\rho$ is called
reducible (irreducible)
if its action on $\C^n$ is reducible (irreducible).  A representation
$\rho$ is called
semi-simple
if it is a direct sum of irreducible representations.
\begin{defin}
$$
\C {\cal M}_B = \{\sigma \in \Hom(\pi_1(X), \GLn) : \sigma
\mbox{ is semi-simple.}\}/\GLn.
$$
$$
{\cal M}_B = \{\sigma \in \Hom(\pi_1(X), \Upo) : \sigma
\mbox{ is semi-simple.}\}/\Upo.
$$
\end{defin}

Let $E$ be a rank $n$ complex
vector bundle over $X$ with 
$\deg(E) = 0$.  Denote by $\Omega$ the 
canonical bundle on $X$.  A holomorphic structure
$\db$ on $E$ induces holomorphic structures on the bundles
$End(E)$ and $End(E) \otimes \Omega$.
A Higgs bundle is a pair $(E_{\db}, \Phi)$,
where $\db$ is a holomorphic structure on $E$ and 
$\Phi \in \H^0(X, End(E_{\db}) \otimes \Omega)$.  Such 
a $\Phi$ is called a Higgs field.
We denote the holomorphic bundle $E_{\db}$ by $V$.

Define the slope of a vector bundle $V$ to be
$$
s(V) = \deg(V)/\rk(V).
$$
For a fixed $\Phi$, a holomorphic sub-bundle
$W \subset V$ is said to be $\Phi$-invariant if
$\Phi(W) \subset W \otimes \Omega$.
A Higgs bundle $(V,\Phi)$ is stable (semi-stable) if $W \subset V$ being
$\Phi$-invariant implies
$$
s(W) < (\le) s(V).
$$
A Higgs bundle is called poly-stable if it is a direct
sum of stable Higgs bundles of the same slope \cite{Hi1,Si2}.

The Dolbeault moduli space
$\C {\cal M}$ is the moduli space 
of poly-stable (or $S$-equivalence 
classes of semi-stable \cite{Ni1}) Higgs bundles
on $X$ \cite{Hi1, Hi2, Ni1, Si2}.  
A Higgs bundle is called reducible if it is poly-stable (semi-stable) 
but not stable.

Now we summarize the relation between the moduli spaces
$\C {\cal M}_{B}$ and $\C {\cal M}$.
\begin{defin} \label{def:M}
Let ${\cal M}$ be the subset of  $\C {\cal M}$ consisting
of Higgs bundles $(V, \Phi)$ satisfying
the following two conditions:
\begin{enumerate}
\item $V$ is a direct sum:
$$
V = V_P \oplus V_Q,
$$
where $V_P, V_Q$ are of ranks $p,1$, respectively.
\item The Higgs field decomposes into two maps:
$$
\Phi_1 : V_P \longrightarrow V_Q \otimes \Omega,
$$
$$
\Phi_2 : V_Q \longrightarrow V_P \otimes \Omega.
$$
\end{enumerate}
\end{defin}
Hence each
$(V,\Phi) = (V_P \oplus V_Q, \Phi) \in {\cal M}$
is associated with an invariant
$$
d = \deg(V_P) = -\deg(V_Q).
$$
The Toledo invariant is defined to be $\tau = 2 d$.
The subset of ${\cal M}$ consisting of classes with a fixed Toledo
invariant $\tau$ is denoted by ${\cal M}_{\tau}$.

\begin{thm} \label{thm:summary}
\mbox{}

\begin{enumerate}
\item The moduli spaces $\C {\cal M}_{B}$ and
$\C {\cal M}$ are homeomorphic.
\item The reducible representations in $\C {\cal M}_{B}$ correspond to
the poly(semi)-stable, but not stable, points.
\item The subspace ${\cal M}_{B}$ is homeomorphic to ${\cal M}$.
\item $\tau \in 2 \Z$ and 
$-2(g-1) \le \tau \le 2(g-1)$.  In particular, equality holds
only if the corresponding Higgs bundle is poly(semi)-stable,
but not stable.
\item The space ${\cal M}_{\tau}$ is homeomorphic to
${\cal M}_{-\tau}$.
\end{enumerate}
\end{thm}
\begin{pf}
See \cite{Co1, Hi1, Si1, Si2, To1, To2, Xi3}.
\end{pf}
By Theorem~\ref{thm:summary} (4) (5), we assume, for the rest of the paper
that 
$$0 \le \tau \le 2g-2.$$

\section{The $\C^*$-Action and the Hodge Bundles}
If $(V,\Phi) \in \C {\cal M}$, then
$(V, t \Phi) \in \C {\cal M}$ for all $t \in \C^*$.
This defines an action \cite{Hi1,Hi2,Si2}
$$
\C^* \times \C {\cal M} \longmapsto \C {\cal M}.
$$
By Definition~\ref{def:M}, we have
\begin{prop} \label{prop:preserve}
The $\C^*$-action preserves ${\cal M}$.
\end{prop}

A Hodge bundle on $X$
is a direct sum of bundles \cite{Si2}
$$
V = \bigoplus_{s,t} V^{s,t}
$$
together with maps (Higgs field)
$$
\Phi_i : V^{s,t} \longrightarrow V^{s-1, t+1} \otimes \Omega.
$$
Definition~\ref{def:M} implies that
\begin{prop} \label{prop:hodge}
Suppose $(V_P \oplus V_Q, (\Phi_1, \Phi_2)) \in  {\cal M}$.
Then $(V_P \oplus V_Q, (\Phi_1, \Phi_2))$ is a Hodge bundle
if and only if $(V_P \oplus V_Q, (\Phi_1, \Phi_2))$
is either binary or ternary in the following sense:
\begin{enumerate}
\item Binary:  $\Phi_2 \equiv 0$.
\item Ternary:  $V_P = V_1 \oplus V_2$ and
the Higgs field consists of two maps:
$$
\Phi_1 : V_2 \longrightarrow V_Q \otimes \Omega,
$$
$$
\Phi_2 : V_Q \longrightarrow V_1 \otimes \Omega.
$$
\end{enumerate}
\end{prop}

\begin{prop} \label{prop:fixed}
A Higgs bundle $(V,\Phi) \in \C {\cal M}$ is a Hodge bundle if and only
if $(V, \Phi) \cong (V, t\Phi)$ for all $t \in \C^*$. 
\end{prop}
\begin{pf}
See \cite{Si2, Si3, Si4}.
\end{pf}

\begin{lemma} \label{lem:extend}
The $\lim_{t \rightarrow 0} (V, t\Phi)$ exists in ${\cal M}$ for any
$(V, \Phi)$ in ${\cal M}$.  In other words, the $\C^*$-action always
extends to a $\C$-action.
\end{lemma}
\begin{pf}
The $\lim_{t \rightarrow 0} (V, t\Phi)$ exists in $\C {\cal M}$
for $(V, \Phi)$ in $\C {\cal M}$ \cite{Si4}.  Since $\Upo$ is a closed
subgroup of $\GLn$,
${\cal M}_B$ is a closed subset of $\C {\cal M}_B$ and the
embedding is proper.  The lemma then follows from Theorem~\ref{thm:summary}.
\end{pf}

From Proposition~\ref{prop:fixed} and Lemma~\ref{lem:extend} and the
facts that the $\C^*$-action preserves $\M$ and $\M$ is closed in
$\C\M$, we have
\begin{cor} \label{cor:deform}
Every Higgs bundle in $\M$ can be deformed to a Hodge bundle in $\M$.
\end{cor}

\begin{lemma} \label{lem:stratum}
Let $E = (V_1 \oplus V_2 \oplus V_Q ,(\Phi_1, \Phi_2))$ be a stable ternary 
bundle in ${\cal M}$ with 
$\Phi_2 \not\equiv 0$.  Then there is a stable Higgs bundle 
$F = (V,\Phi) \in {\cal M}$ not isomorphic to $E$ such
that 
$$
\lim_{t \rightarrow \infty} (V, t \Phi) = E.
$$
\end{lemma}
\begin{pf}
This is essentially Lemma 11.9 in \cite{Si4}.  There is one additional 
and crucial observation to be made.  That is the Higgs bundle $F$ in 
the proof is actually in ${\cal M}$.  This is immediate from the 
construction of $F$.
\end{pf}

\begin{prop} \label{prop:binary}
Every Higgs bundle in ${\cal M}$ can be deformed to a binary Hodge bundle.
\end{prop}
\begin{pf}
The proof parallels the proof of Corollary 11.10 in \cite{Si4}.
The only difference is that for the subset ${\cal M}_{\tau}$, the lowest
stratum is the space of binary Hodge bundles with Toledo invariant
$\tau$, instead of $M_{p+1,0}$.  Suppose $(V,\Phi) \in \M_{\tau}$.
Then it is of the form
$(W',0) \oplus (W',\Phi')$ where $(W',\Phi')$ is stable.  
By Corollary~\ref{cor:deform}, we may assume $(W',\Phi')$ is
ternary.  Then Lemma~\ref{lem:stratum} implies that
$(W',\Phi')$ is not in the lowest
possible stratum, hence, can be deformed to a 
fixed point set (with respect to the $\C^*$-action) of the lowest stratum
which consists of binary Higgs bundles (See the section titled
``Actions of $\C^*$'' in \cite{Si4} for further discussions of these strata).  
\end{pf}
\begin{defin}
Let $B_{\tau}$ be the space of all poly-stable (or $S$-equivalence classes
of semi-stable \cite{Ni1}) binary Hodge bundles
$(V_P \oplus V_Q,(\Phi_1,0))$
with $\deg(V_P) = d = - \deg(V_Q)$ and 
$\tau = 2 d$.
\end{defin}
The rest of the paper is devoted to showing that $B_{\tau}$ is 
connected.

\section{The Deformation of Binary Hodge Bundles}
A family (or flat family) of Higgs bundles $(V_Y, \Phi_Y)$ 
is a variety Y such that there is a vector bundle $V_Y$ on 
$X \times Y$ together with a section
$\Phi_Y \in \Gamma(Y, (\pi_Y)_*(\pi_X^* \Omega \otimes End(V_Y))$ \cite{Ni1}.
$\C {\cal M}$ being a moduli space implies that
if $Y$ is a family of stable (poly-stable or 
$S$-equivalence classes of semi-stable) 
Higgs bundles, then there is a natural morphism \cite{Mu1, Ne1}
$$
t : Y \longrightarrow \C {\cal M}.
$$
Moreover $t$ takes every point $y \in Y$ to the point of
$\C {\cal M}$ that corresponds to the Higgs bundle
in the family over $y$ \cite{Mu1, Ne1, Ni1}.

The space ${\cal M}$ is a subvariety of 
$\C {\cal M}$;
hence, to show that two stable (poly-stable or 
$S$-equivalence classes of semi-stable) Higgs bundles 
$(V_1,\Phi_1)$ and $(V_2,\Phi_2)$
belong to the same component of ${\cal M}$, it suffices to
exhibit a connected family $Y$
(whose image under the morphism $t : Y \longrightarrow \C {\cal M}$ is
contained in ${\cal M}$) of stable
(poly-stable or $S$-equivalence classes of semi-stable) 
Higgs bundles containing both 
$(V_1,\Phi_1)$ and $(V_2,\Phi_2)$. 

Let $(V_P \oplus V_Q, \Phi) \in B_{\tau}$ be a stable 
binary Hodge bundle with
$d = \deg(V_P) = - \deg(V_Q)$.

\subsection{The Grothendieck Quot scheme}
Denote by $H_{r,d_1}$ the set of all vector bundles
of rank $r$ and degree $d_1 \le 0$ with the property that
if $W \in H_{r,d_1}$ and $U \subset W$, then
$\deg(U) \le 0$.
\begin{prop} \label{prop:quot}
Suppose $W \in H_{r,d_1}$.
Then for any line bundle $L$ with $\deg(L) > 2g-1 - d_1$,
\begin{enumerate}
\item $\H^1(W \otimes L) = 0$,
\item $W \otimes L$ is generated by global sections.
\end{enumerate}
\end{prop}
\begin{pf}
The proof resembles the proof of Lemma 20 in Chapter I of \cite{Se1}.
By Serre duality, 
$\H^1(W \otimes L) \cong \H^0(\Omega \otimes W^* \otimes L^*)$.
Hence if $\H^1(W \otimes L) \neq 0$, then
$\Omega \otimes W^* \otimes L^*$ contains a line bundle $N$ of
degree greater than or equal to 0: 
$$
0 \longrightarrow N \longrightarrow \Omega \otimes W^* \otimes L^*
\longrightarrow (\Omega \otimes W^* \otimes L^*)/N \longrightarrow 0.
$$
Dualizing and tensoring with $\Omega \otimes L^*$, 
we obtain an exact sequence
$$
0 \longrightarrow U \longrightarrow W
\longrightarrow L' \longrightarrow 0,
$$
with $\deg(U) \ge d_1 - (2g-2) + \deg(L)$.
Since $W \in H_{r,d_1}$ and $U \subset W$, $\deg(U) \le 0$.  This implies
$\deg(L) \le 2g-2 - d_1$ which is a contradiction.  This shows that
(1) is true for any $L > 2g - 2 - d_1$.

The proof of (2) essentially reduces to
showing that 
$$
\H^1(W \otimes L \otimes L_x^{-1}) = 0
$$
(see the proof of Lemma 20 in Chapter I of \cite{Se1}), where
$L_x$ is the ideal sheaf at a point $x \in X$.  Since $\deg(L_x) = 1$,
$\deg(L \otimes L_x^{-1}) > 2g-2 - d_1$.  Since (1) is true for any
$L > 2g - 2 - d_1$,  (2) follows.
\end{pf}
Let $D = 2gr + (1-r)d_1$
and $a = D + r(1-g) = r(g+1) + (1-r)d_1$.  For the pair $(a,r)$, we construct
the Grothendieck scheme $Q$ parameterizing
the quotient sheaves of $\O^a$ with Hilbert polynomial
$H(m) = a + rm$ \cite{Gr1}.  The scheme $Q$ 
contains the sub-scheme $R$ defined by:
$$
R = \{W \in Q : W \mbox{ is locally free and } \H^1(W) = 0 \}.
$$
The sub-scheme $R$ is smooth and connected (The proof is
the same as that of Proposition 23 in Chapter I of \cite{Se1}).

Suppose $L$ is a line bundle of degree $-2g+d_1$.  If $W \in R$,
then $\deg(W \otimes L) = d_1$. 
Hence $R$ also parameterizes a family of 
vector bundles of degree $d_1$ and rank $r$.  By Proposition~\ref{prop:quot},
$R$ contains all the bundles in $H_{r,d_1}$.
We shall denote the scheme $R$ so constructed as $R_{r,d_1}$.

\subsection{The canonical factorization}  Let 
$(V_P \oplus V_Q, \Phi) \in B_{\tau}$ be a binary Higgs bundle 
with $\Phi \not\equiv 0$.
There exist
bundles $V_1, V_2$ and $W_1, W_2$ such that
the following diagram (the canonical factorization \cite{Sh1})
$$
\begin{array}{ccccccccc}
0 & \longrightarrow & V_1 & \stackrel{f_1}{\longrightarrow} &
V_P & \stackrel{f_2}{\longrightarrow} &
V_2 & \longrightarrow & 0\\
  & & & & \Phi \Big\downarrow & &
\varphi \Big\downarrow & & \\
0 & \longleftarrow & W_2 & \stackrel{g_2}{\longleftarrow}
& V_Q \otimes \Omega & \stackrel{g_1}{\longleftarrow} &
W_1 & \longleftarrow & 0
\end{array}
$$
commutes, and the rows are exact, $\rk(V_2) = \rk(W_1)$ and
$\varphi$ has full rank at a generic point of $X$.
As $V_Q \otimes \Omega$ is a line bundle, $V_Q \otimes \Omega = W_1$.
Thus the diagram simplifies to
$$
\begin{array}{ccccccccc}
0 & \longrightarrow & V_1 & \stackrel{f_1}{\longrightarrow} &
V_P & \stackrel{f_2}{\longrightarrow} &
V_2 & \longrightarrow & 0\\
  & & & & \Phi \Big\downarrow & &
\varphi \Big\downarrow & & \\
& & &
& V_Q \otimes \Omega & = &
W_1 & &
\end{array}
$$
Let $d_1 = \deg(V_1)$ and $d_2=\deg(V_2)$.
Since $(V_P \oplus V_Q, \Phi)$ is semi-stable,
$d_1/(p-1) \le 0$.
Since $\varphi \not\equiv 0$, $d_2 \le (2g-2) - d$ (Recall that
$d = \deg(V_P)$).
Hence $d_1 \ge 2d - (2g-2)$.
To summarize
$$
\left\{
\begin{array}{lllll}
2d - (2g-2) & \le & d_1 & \le & 0\\[2ex]
d & \le & d_2 & \le & -d + (2g-2)
\end{array}
\right.
$$
Denote by $B_{\tau}(d_2)$ the subspace of $B_{\tau}$ such that
$(V, \Phi) \in B_{\tau}(d_2)$ implies $\deg(V_2) = d_2$ in the above
canonical factorization.
\begin{prop} \label{prop:strata}
The space $B_{\tau}(d_2)$ is connected.
\end{prop}
\begin{pf}
Denote by $J^l$ the Jacobi variety identified with
the set of holomorphic line bundles of degree $l$ on $X$.
For each $V_2 \in J^{d_2}$, the variety
$\C^* \times X^{-d + 2(g-1) - d_2}$
parameterizes a family of pairs that contains
all pairs $(V_Q, \varphi)$ such that
$V_Q \in J^{-d}$ and
$$
0 \not\equiv \varphi \in \H^0(X, V_2^{-1} \otimes V_Q \otimes \Omega).
$$
Note the moduli of all such pairs is simply 
$\C^* \times Sym^{-d + 2(g-1) - d_2} X$,
where $Sym^t X$ is
the $t$-th symmetric product of $X$, i.e. the moduli is 
$\C^* \times X^{-d + 2(g-1) - d_2}$ quotiented by the
symmetry group on $-d + 2(g-1) - d_2$ letters.
Hence the variety
$$
S = J^d \times (\C^* \times X^{-d + 2(g-1) - d_2})
$$
parameterizes a family of triples that contains all
triples $(V_2, V_Q, \varphi)$ such that
$$
V_2 \stackrel{\varphi}{\longmapsto} V_Q \otimes \Omega,
$$
with $V_2 \in J^{d_2}, V_Q \in J^{-d}$ and $\varphi \not\equiv 0$.
The variety $S$ is smooth.

Suppose $d_2 = 0$.  Then $d_1 = 0$ and every Higgs bundle in $B_{\tau}(0)$
is reducible.  There are two cases.

Case 1: $\tau > 0$.  Then $\Phi \not\equiv 0$ and 
every Higgs bundle in $B_{\tau}(0)$ is
of the form $(V_1,0) \oplus (V_2 \oplus V_Q, \Phi)$, 
where $V_1 \in M_{p-1,0}, V_2 \in J^d, V_Q \in J^{-d}$.
Hence $M_{p-1,0} \times \C^* \times Sym^{-2d + (2g-2)}X \times J^{-d}$
parameterizes a family that contains
every Higgs bundle in $B_{\tau}(0)$.  Since 
$M_{p-1,0} \times \C^* \times Sym^{-2d + (2g-2)}X \times J^{-d}$
is connected, $B_{\tau}(0)$ is connected.

Case 2:  $\tau = 0$.
Then every Higgs bundle in $B_0(0)$ must be one of the following two forms:
\begin{enumerate}
\item $(V_P \oplus V_Q, 0)$, where $V_P \in M_{p,0}, V_Q \in J^0$,
\item $(V_1,0) \oplus (V_2 \oplus V_Q, \Phi)$, where 
$V_1 \in M_{p-1,0}, V_2 \in J^0, V_Q \in J^0$. 
\end{enumerate}
Type 2 can be deformed to type 1 by simply deforming the
Higgs field $\Phi$ to zero.  Since $M_{p,0} \times J^0$
is connected, $B_0(0)$ is connected.

\begin{lemma} \label{lem:constant}
Suppose $d_2 > 0$.
Then the dimension of the space $\Ext^1(V_2, V_1)$ is 
$(p-1)(g-1) + (p-1)d_2 - d_1$.
\end{lemma}
\begin{pf}
The subspace $V_1$ is $\Phi$-invariant.  By semi-stability,
$V' \subset V_1$ implies
$$
s(V') \le 0 < d_2 = \deg(V_2).
$$
Hence $\H^0(\Hom(V_2, V_1)) = 0$.
The lemma then follows from the fact that 
$\Ext^1(V_2,V_1) \cong \H^1(\Hom(V_2,V_1))$ and from 
Riemann-Roch.
\end{pf}
For $d_2 > 0$, construct the
universal bundle \cite{At1, Se1}
$$
U \longrightarrow X \times R_{p-1,d_1} \times S
$$
such that 
$$
U|_{(X, V_1, V_2, V_Q, \varphi)} = V_2^{-1} \otimes V_1.
$$ 
Let $\pi$ be the projection
$$
\pi : X \times  R_{p-1,d_1} \times S \longrightarrow R_{p-1,d_1} \times S.
$$
Applying the right
derived functor $R^1$ to $\pi$ gives
the sheaf $\F = R^1\pi_*(U)$ \cite{Ha1} such that
$$
\F|_{(V_1,V_2,V_Q,\varphi)} = \H^1(X, V_2^{-1} \otimes V_1).
$$
By Lemma~\ref{lem:constant} and Grauert's theorem \cite{Ha1},
$\F$ is locally free, hence, is associated with a vector bundle
$$
F \longmapsto R_{p-1,d_1} \times S
$$
of rank $(p-1)(g-1) + (p-1)d_2 - d_1$.  Since $R_{p-1,d_1}$ is
smooth and connected, the total
space $F$ is smooth, connected and parameterizes a family of Higgs bundles 
that fit into the canonical decomposition with fixed $d_2$.
By construction, the scheme $F$
parameterizes a family of Higgs bundles that contains
every member in the parameter space $B_{\tau}(d_2)$.  Moreover
if a Higgs bundle in $F$ is semi-stable, it must belong 
to $B_{\tau}(d_2)$.
Since the semi-stability condition is open \cite{Si3},
the subset of $F$ parameterizing the semi-stable 
Higgs bundles, if not empty, is open and dense in $F$, hence, 
connected.  This implies $B_{\tau}(d_2)$ is connected.
\end{pf}

\subsection{Deformation between the $B_{\tau}(d_2)$'s}
Fix a set of distinct points 
$$
A = \{x_1,...,x_{d_2}, y_1,...,y_{d_2-1}, z_1,...,z_{d_2-d-1}\} \subset X
$$
and let $Y = X \setminus A$.  Fix $y \in Y$.
For $t \in Y$, consider the following divisors on $X$:
\begin{eqnarray*}
D_2 & = & \{x_1,...,x_{d_2}\},\\
C(t) & = & \{x_1,...,x_{d_2}, -y_1,...,-y_{d_2-1},-t \},\\
C & = & \{x_1,...,x_{d_2},-z_1,...,-z_{d_2-d-1}, -y \}.
\end{eqnarray*}
For any divisor $D$, denote by $[D]$ the line bundle associated with
$D$.
The set $Y$ parameterizes a family of Higgs bundles as follows.  Let
$$
V_P(t) = [C] \bigoplus_{i=1}^{p-1} [C_i(t)],
$$
where $C_i(t) = C(t)$ for all $i$ and denote the projection maps to
the $[C]$ and $[C_i(t)]$ factors by
$p$ and $p_i(t)$, respectively.  The divisors
$D_2 - C(t)$ and $D_2 - C$ define maps
$h_i(t): [C_i(t)] \longrightarrow [D_2]$ and 
$h : [C] \longrightarrow [D_2]$, respectively.  These maps induce a map
$$
G_t : V_P(t) \longrightarrow [D_2], \ \ G_t = h + \sum_{i=1}^{p-1} h_i(t).
$$
Let $V_2 = [D_2]$.  Since $d_2 \le (2g - 2) - d$, there exists 
$V_Q \in J^{-d}$
and $0 \not\equiv \varphi \in \H^0(V_2^{-1} \otimes V_Q \otimes \Omega)$.
Let $\Phi(t) = \varphi \circ G_t$.  Then
$(V_P(t) \oplus V_Q, \Phi(t))$ is a family of Higgs bundles parameterized
by $Y$.  Let $p_P, p_Q$ be the projections onto 
the $V_P(t),V_Q$ factors.

\begin{lemma} \label{lem:subsheaf}
If $U \subset V_P(t)$, then $\deg(U) \le d$.
\end{lemma}
\begin{pf}
This is an inductive argument.  

Case 1: $p(U) \neq 0$.  Consider
the sequence
$$
0 \longrightarrow U' \longrightarrow U \longrightarrow p(U) \longrightarrow
0.
$$
Then $\deg(p(U)) \le d$ and $\deg(U) \le \deg(U') + d$.
Now we begin with the smallest $i$ with $p_i(t)(U') \neq 0$ and construct
$$
0 \longrightarrow U_i \longrightarrow U' \longrightarrow p_i(U') \longrightarrow
0.
$$
Again $\deg(p_i(U')) \le \deg([C_i(t)]) = 0$ and $\deg(U) \le \deg(U_i) + d$.
Now we let $j>i$ be the smallest integer with $p_j(U_i) \neq 0$ and
construct the new sequence and obtaining $U_j$ with
$\deg(U) \le \deg(U_j) + d$.  Note $\rk(U) > \rk(U') > \rk(U_i)$,
so eventually the process ends and since $\deg([C_i(t)]) = 0$ for all $i$,
we have $\deg(U) \le d$.

Case 2: $p(U) = 0$.   Here we simply begin with the smallest $i$ with 
$p_i(t)(U') \neq 0$ as in Case 1.  The rest is the same and we conclude
that $\deg(U) \le 0 \le d$.
\end{pf}
\begin{lemma} \label{lem:small}
Suppose $L \subset V_P(t)$ is a line bundle with $\deg(L) > 0$,
then $L = [C]$.
\end{lemma}
\begin{pf}
Suppose $p_i(t)(L) \not\equiv 0$ for some $i$, then 
$\deg(L) \le \deg(C_i(t)) = 0$, a contradiction.
\end{pf}
\begin{prop} \label{prop:deform}
The Higgs bundle $(V_P(t) \oplus V_Q, \Phi(t))$ is in $B_{\tau}(d_2-1)$ if
$t = y$ and in $B_{\tau}(d_2)$ if $t \neq y$. 
\end{prop}
\begin{pf}
From the definition, one only needs to check that 
$(V_P(t) \oplus V_Q, \Phi(t))$ is semi-stable for all $t \in Y$.
Suppose $W \in V_P(t) \oplus V_Q$ is $\Phi(t)$ invariant.
There are two cases.

Case 1: $p_Q(W) \neq 0$.  Then there is an exact sequence
$$
0 \longrightarrow U \longrightarrow W \longrightarrow p_Q(W) \longrightarrow
0.
$$
Since $U \subset V_P$, by Lemma~\ref{lem:subsheaf}, 
$\deg(U) \le d$.  Since $\deg(p_Q(W)) \le \deg(V_Q) = -d$, 
$\deg(W) \le d - d = 0$.

Case 2: $p_Q(W) = 0$.  Since $W$ is $\Phi(t)$-invariant,
$W \subset \ker(G_t)$.
It is immediate from the definition that, $[C] \not\subset \ker(G_t)$.
Now we begin the construction similar to that in the proof of 
Lemma~\ref{lem:subsheaf}.
Begin with the smallest $i$ with $p_i(t)(W) \neq 0$ and construct
$$
0 \longrightarrow U_i \longrightarrow W \longrightarrow p_i(W) \longrightarrow
0.
$$
As $\deg(p_i(W)) \le \deg([C_i(t)]) = 0$, $\deg(W) \le \deg(U_i)$.
Continue this process as before.  
Eventually we reach the exact sequence
$$
0 \longrightarrow U_k \longrightarrow U_j \longrightarrow p_k(U_j) 
\longrightarrow 0,
$$
with 
$$\deg(W) \le \deg(U_j) = \deg(U_k) + \deg(p_k(U_j)) \le \deg(U_k)$$ 
and $U_k$ a line bundle.
Since $[C] \not\subset \ker(G_t)$ and $U_k \subset \ker(G_t)$, 
by Lemma~\ref{lem:small}, $\deg(U_k) \le 0$.
Hence $\deg(W) \le \deg(U_k) \le 0$.
Hence we conclude that the Higgs bundle $(V_P(t) \oplus V_Q, \Phi(t))$
is semi-stable.
\end{pf}
Proposition~\ref{prop:binary} states that every Higgs bundle in
${\mathcal M}$ may be deformed to a Higgs bundle in $B_{\tau}$ which
is the union of the $B_{\tau}(d_2)$'s.  Proposition~\ref{prop:strata}
asserts that each $B_{\tau}(d_2)$ is connected.  Theorem~\ref{thm:main}
then follows from Proposition~\ref{prop:deform}.

\end{document}